\documentclass{amsart}
\usepackage{amscd,amssymb}

\theoremstyle{plain}
\newtheorem{thm}{Theorem}[section]
\newtheorem{lem}[thm]{Lemma}

\newtheorem{prop}[thm]{Proposition}

\theoremstyle{definition}
\newtheorem{Def}[thm]{Definition}

\theoremstyle{remark}

\newtheorem{Remark}[thm]{Remark}
\newtheorem{Assumption}[thm]{Assumption}

\numberwithin{equation}{section}

\begin{document} 

\textheight = 650pt

\setcounter{section}{-1}

\title{Non-nesting actions of groups on real trees vs isometric actions}
\author{A.Ivanov}
\thanks{A.Ivanov, \\ 
Institute of Mathematics, University of Wroc{\l}aw, \\ 
pl.Grunwaldzki 2/4, 50-384 Wroc{\l}aw, Poland, \\ 
e-mail : ivanov@math.uni.wroc.pl \\ 
Fax number: +48-71-3757429  }

\maketitle

\bigskip
{\bf Abstract} 
We study non-nesting actions of groups on
$\mathbb{R}$-trees.
We prove that some natural conditions describing  
how the group is generated, imply that such an action 
involves an isometric action on an $\mathbb{R}$-tree. 
This can be applied to permutation groups, linear groups 
and some Polish groups.\parskip0pt



\section{Introduction}

Non-nesting actions by homeomorphisms on $\mathbb{R}$-trees
frequently arise in geometric group theory 
(for example see \cite{bow} and \cite{DrSa}).
Explicitly they were introduced in \cite{levitt}.
Following that paper, we concentrate on the question when 
a non-nesting action of a group on an $\mathbb{R}$-tree 
can be 'transformed' into an isometric action. 

We introduce an algebraic condition on a group 
which in the situation when $G$ acts on a real tree $T$ reduces 
this question to {\em characteristic subsets} $C_g\subseteq T$ of 
some special elements $g\in G$. 

\begin{Def} \label{path} 
Let $G$ be a group and $X\subset G$ be a conjugacy class of $G$. \\ 
(1) We say that a sequence $g_1 ,...,g_n$ of elements of $X$ 
is an $X$-{\bf path} if for each pair $g_i$, $g_{i+1}$ 
with $i<n$ there are $\varepsilon$ and $\tau\in \{ -1 ,+1\}$ 
such that $g^{\varepsilon}_{i}\cdot g^{\tau}_{i+1}$ is conjugate  
with some element from $X \cup X^{-1}$. \\ 
(2) We say that the set $X$ is {\bf path connected} 
if for any $g,g'\in X$ 
there is an $X$-path, where the first and the last elements are $g$ and $g'$. 
\end{Def} 

We will show that when $G=\langle X\rangle$ and 
$X$ is path connected, if we choose $v\in T$ in an appropriate way, 
then the the set of all elements of the orbit $Gv\subset T$ 
together with medians of all triples of $Gv$ (denote it by $c(Gv)$), 
can be equivariantly embedded into a $G$-space of an isometric action 
on a real tree. 
This can be applied to Polish groups with a comeagre conjugacy class, 
to algebraically closed groups, to $SL(n,K)$ for various rings $K$, 
to infinite dimensional linear groups, to groups of bounded automorphisms 
of relatively free algebras and to the automorphism group of a regular tree.   

To explain what "an appropriate choice of $v\in T$" means 
note that if $c(Gv)$ is discrete then it can be always extended 
to an isometric action of $G$ on a real tree. 
Otherwise it is possible that the $G$-space $c(Gv)$ is "too poor" 
as in the following situation:  
\begin{quote} 
\begin{itemize} 
\item $G=\langle X\rangle$ and $X$ is path connected;  
\item all elements of $X$ act as loxodromic elements on $T$;   
\item but for any $g\in X$ any orbit in $T$ of the setwise stabiliser 
$G_{\{ L_g\}}$ of the axis of $L_g$ is not dense in $c(Gv)\cap L_g$. 
\end{itemize} 
\end{quote}  
In this case the choice of $v$ is not appropriate, 
because the traces of axises 
in the $G$-space $c(Gv)$ are too small to reconstruct the action. 
In Section 3 we prove a theorem (Theorem \ref{main}) on embedding of $c(Gv)$ 
into $G$-spaces of isometric actions, where we give general conditions 
which make the choice of $v\in T$ acceptable for us.

In Sections 1 and 2 we do some preliminary work. 
In particular since we apply several tools from \cite{bow} used 
for pretrees, in Section 1 we give some introduction to pretrees and flows, 
and prove an important lemma applied in Section 3. 
Section 2 is devoted to non-nesting actions on real trees.  
Here we concentrate on stabilisers of ends. 

The results of the paper are slightly connected with \cite{guiva} 
where actions of Polish groups on trees were studied. 
It is worth noting that some related problems 
(for example of embedding of generalized trees into 
$\mathbb{R}$-trees and $\Lambda$-trees) have been studied before 
(see \cite{bow}, \cite{bc}, \cite{chi}, \cite{dun}, 
\cite{levitt}, \cite{mactho} ).
Our motivation is partially based on these investigations.

\bigskip 

{\bf Applications.} 

We now give a number of examples. where the condition of 
path connectedness naturally arises. 
Thus in all cases below the conclusion of Theorem \ref{main} 
holds when we have an appropriate action. 

\subsection{Polish groups with a comeagre conjugacy class} 
A Polish group is a topological group whose topology is
{\em Polish} (a Polish space is a separable completely
metrizable topological space).
A subset of a Polish space is {\em comeagre} if it containes
an intersection of a countable family of dense open sets.
The group $Sym(\omega )$ of all permutations of a countably infinite set 
is a typical example of a Polish group with a comeagre conjugacy class 
(see also \cite{iva} for some descriptions of permutations of this class). 

Let $X$ be a comeagre conjugacy class of a Polish group $G$. 
Since for every $g\in G$ the intersection $X^{-1} \cap X\cap gX$ 
is also comeagre we see $G=X\cdot X$. 
Moreover since for any two $h_1 ,h_2 \in X$ the intersection 
$h_1 X \cap h_2 X \cap X$ is comeagre we see that $h_1$ and $h_2$ 
are connected by an $X$-path of length 3. 

It is worth noting that Theorem \ref{main} is not appropriate 
in this case. 
By \cite{guiva} if a Polish group $G$ has a comeagre conjugacy class 
then for any non-nesting action on a real tree any element of $G$ 
fixes a point. 
Thus the condition of ${\bf e}$-contractibility of Theorem \ref{main} 
is not satisfied. 
We give this example only for illustration of path connectedness. 

\bigskip 

\subsection{Algebraically closed groups.} 
Let $G$ be a non-trivial group. 
It is called {\em algebraically closed} if any finite system of equations 
with parameters having a solution in some supergroup already 
has a solution in $G$. 
If in this definition we consider finite systems of equations and inequations, 
we obtain the same notion. 
If $G$ is algebraically closed then by Lemma 1 of \cite{mac} the set $X$ of 
all elements of infinite order forms a conjugacy class and by Corollary 2 
of \cite{mac}, $G=X\cdot X$. 

For any $g_1 ,g_2 \in X$ embed the group $G*\langle h \rangle$
with $|h|=\infty$, into another algebraically closed group $G_1$ 
and find $u ,v, w \in G_1$ so that $g_1 h=u^{-1} g_1 u$, 
$h g_2 =v^{-1} g_1 v$ and $h=w^{-1} g_1 w$. 
Since $G$ is algebraically closed there are $h',u',v'$ and $w'$ in $G$ 
so that $g_1 h'=(u')^{-1} g_1 u'$, $h' g_2 =(v')^{-1} g_1 v'$ and 
$h'=(w')^{-1} g_1 w'$.
We see that $g_1 ,h',g_2 $ is an $X$-path of length 3.  

\bigskip 

\subsection{The group $SL(n,K)$, $n>2$, over a field $K$.} 
The algebraic material below is standard. 
We use \cite{vast}. 
For any ring $R$ the group $SL(n,R)$ is generated 
by the family of all elementary transvections 
$$ 
t_{ij}(\alpha ) = I+\alpha I_{ij} \mbox{ , where } \alpha\in K \mbox{ and } 1\le i,j \le n. 
$$ 
A {\em transvection} is a matrix of the form 
$$ 
t_{\bar{u}\bar{v}}(\xi )=I+\bar{u}\xi \bar{v} \mbox{ , where } \bar{u} \mbox{ and } \bar{v} 
\mbox{ are a column and a row so that } \bar{v}\bar{u} = 0. 
$$ 
Since 
$a\cdot t_{\bar{u}\bar{v}}(\xi )\cdot a^{-1} = t_{a\bar{u},\bar{v}a^{-1}}(\xi )$,  
the set of transvections is conjugacy invariant. 
The following equality is a version of the well-known Chevalley commutator formula: 
$$ 
[t_{\bar{u}\bar{v}}(\xi ), t_{\bar{w}\bar{y}}(\zeta )] = t_{\bar{u}\bar{y}}(\xi \bar{v}\bar{w} \zeta )  
\mbox{ , where } \bar{v}\bar{u}=\bar{y}\bar{w}=\bar{y}\bar{u}=0 \mbox{ and } \xi ,\zeta \in R . 
$$ 
Note that 
$t_{\bar{u}\bar{v}}(\xi ) = t_{\xi\bar{u},\bar{v}}(1 )= t_{\bar{u},\xi\bar{v}}(1 )$ 
and 
$t_{\bar{u}+\bar{w},\bar{v}}(\xi ) = t_{\bar{u}\bar{v}}(\xi )\cdot t_{\bar{w}\bar{v}}(\xi )$,
where $\bar{v}\bar{u} = \bar{v}\bar{w}=0$. 

When $R=K$ is a field, all transvections belong to the same conjugacy class. 
Let $X$ be this class. 
Then $\langle X\rangle = G$. 

{\bf Claim.} {\em Any two transvections of $X$ are connected by an $X$-path of length 5.} 

Indeed, given $t_{\bar{u}\bar{v}}(\xi )$ and $t_{\bar{u}'\bar{v}'}(\xi ')$ 
we find a column $\bar{x}$ and a row $\bar{y}$ with 
$\bar{y}\bar{u} =\bar{v}'\bar{x} = \bar{y}\bar{x} =0$. 
From the Chevalley commutator formula we see   
$$ 
t_{\bar{u}\bar{v}}(\xi )\cdot t_{\bar{u}\bar{y}}(\xi \bar{v}\bar{x} ) = t_{\bar{x}\bar{y}}(-1 ) t_{\bar{u}\bar{v}}(\xi ) t_{\bar{x}\bar{y}}(1 )  
\mbox{ and }
$$ 
$$
t_{\bar{x}\bar{y}}(1)\cdot t_{\bar{x}\bar{v}'}(\xi' \bar{y}\bar{u}' ) = t^{-1}_{\bar{u}'\bar{v}'}(\xi' ) t_{\bar{x}\bar{y}}(1 ) t_{\bar{u}'\bar{v}'}(\xi' )  
$$ 
Thus we can connect $t_{\bar{u}\bar{v}}(\xi )$ and $t_{\bar{u}'\bar{v}'}(\xi ')$ by the following $X$-path  
$t_{\bar{u}\bar{v}}(\xi )$, $t_{\bar{u}\bar{y}}(\xi \bar{v}\bar{x})$, $t_{\bar{x}\bar{y}}(1 )$,  
$t_{\bar{x}\bar{v}'}(\xi'\bar{y}\bar{u}' )$ and $t_{\bar{u}'\bar{v}'}(\xi ')$. 
It is worth noting here that both 
$t_{\bar{u}\bar{y}}(\xi \bar{v}\bar{x})\cdot t_{\bar{x}\bar{y}}(1 )$ and 
$t_{\bar{x}\bar{v}'}(\xi'\bar{y}\bar{u}' ) \cdot t_{\bar{u}'\bar{v}'}(\xi' )$ 
are transvections by appropriate versions of formulas above. 

\bigskip

\subsection{Infinite dimensional general linear groups.} 
Let $V$ be an infinite-dimensional vector space over a division ring $D$. 
A subspace $U$ is called {\em moietous} if $dim U= dim V = codim U$. 
An automorphism $\pi \in GL(V)$ is called a {\em moietous involution} 
if there is a decomposition $V=U\oplus U'\oplus W$ into a sum of moietous 
subspaces so that $\pi$ is identity on $W$ and is an involution 
permuting $U$ and $U'$. 
It is proved in \cite{tolstych1} that $GL(V)$ is generated by the 
conjugacy class of $\pi$. 
Denote this class by $X$.  
We state that 
{\em any two moietous involutions are connected by an $X$-path of length 5}. 

To see this statement take two moietous involutions $\pi_1$ and $\pi_2$ 
with the corresponding decompositions $V =U_1 \oplus \pi_1 (U_1 )\oplus W_1$ 
and $V =U_2 \oplus \pi_2 (U_2 )\oplus W_2$. 
Taking a moietous subspace $U'_1$ of $U_1$ we obtain a decomposition 
$V= U'_1 \oplus \pi_1 (U'_1 ) \oplus W'_1$ corresponding to 
the involution $\pi'_1$ which coincides with $\pi_1$ on $U'_1$ 
and is identical on $W'_1$. 
It is clear that the product $\pi_1 \cdot \pi'_1$ is a 
moietous involution too. 

Taking appropriate $\pi'_1$ and $\pi'_2$ if necessary, 
we may assume that $dim(W_1 \cap W_2 )= dim V$.  
Decomposing $W_1 \cap W_2= U''_1 \oplus U''_2\oplus U''_3$ into a sum 
of moietous subspaces define a moietous involution $\pi_3$ 
which permutes $U''_1$ and $U''_2$ and is identical on the 
remaining part of $V$. 
It is clear that $\pi_1 \cdot \pi_3$ and $\pi_3 \cdot \pi_2$ 
are moietous involutions. 
\bigskip 

\subsection{Subgroups of bounded automorphisms of relatively free algebras.} 
Let $\mathcal{C}$ be a variety of algebras and let 
$F$ be a free $\mathcal{C}$-algebra of uncountable rank $\varkappa$. 
Take a basis $B$ of $F$ and an infinite cardinal $\kappa <\varkappa$. 
For every $\alpha \in Aut(F)$ the $B$-{\em support} of $\alpha$  
is the set $supp(\alpha )= \{ b\in B: \alpha (b)\not= b\}$. 
Define 
$$
Aut_{\kappa ,B}(F) = \{ \alpha \in Aut(F) : supp (\alpha )\le \kappa \}. 
$$  
Let $\pi\in Aut_{\kappa ,B}(F)$ be an automorphism induced by 
a permutation of some basis $B_1$. 
If there is $B_0 \subset B_1$ of cardinality $\kappa$ so that 
$\pi$ acts on $B_0$ by an involution without fixed points and 
fixes $B_1 \setminus B_0$ pointwise, then we say that 
$\pi$ is a {\em moietous involution}.    
It is proved in Theorem 1.5 of \cite{tolstych2} that $Aut_{\kappa ,B}(F)$ 
is generated by the conjugacy class $X$ of moietous involutions. 

The argument of the previous example can be modified to prove that any two 
moietous involutions are connected by an $X$-path of length $\le 5$. 
\bigskip 

\subsection{The automorphism group of a regular tree.} 
Let $A$ be a regular simplicial tree and the valency 
$k$ of a vertex of $A$ is not less than $4$. 
We decompose $A$ as follows. 
Let $e_0 =(a_0 ,a_1 )$ be an edge of $A$ and $T_0$ be 
a rooted subtree of $A$ consisting 
of all elements $v\in A$ so that 
the shortest path from $v$ to $a_0$ does not contain $a_1$. 
It is clear that $a_0$ is the root of $T_0$ and $a_1$ is 
the root of the subtree $A\setminus T_0$.  
Moreover all elements of $A\setminus (T_0 \cup \{ a_1 \} )$ are 
of valency $k$ in $A\setminus T_0$. 
Let $b_1 ,...,b_{k-1}$ be all neighbours of $a_1$ distinct from 
$a_0$ and let $T_1 ,...,T_{k-1}$ be the corresponding rooted 
subtrees of $A\setminus T_0$ with roots $b_i$ (without $a_1$). 
We consider $A$ as $\{ a_1 \}\cup \bigcup_{0\le i\le k-1} T_i$. 

We fix isomorphisms $\phi_i : T_0 \rightarrow T_i$ taking $a_0$ to $b_i$. 
For each pair $(i,j)$ let $\phi_{i,j} = \phi_j \cdot \phi^{-1}_i$. 
It is clear that for every triple $i,j,k$, 
$\phi_{i,j} = \phi_{k,j} \cdot \phi_{i,k}$. 

We also fix an isomorphism $\tau: A\setminus T_0 \rightarrow T_0$ 
taking $a_1$ to $a_0$. 
Let $T'_1 ,...,T'_{k-1}$ be subtrees of $T_0$ corresponding 
to $T_1 ,...,T_{k-1}$ by $\tau$ and let $b'_1 ,...,b'_{k-1}$ 
be the corresponding roots. 
Let $\alpha$ be an automorphism of $T_0$ of order 3 
which fixes $T'_4 \cup ...\cup T'_{k-1}$ pointwise and defines 
the permutation $T'_1 \rightarrow T'_2 \rightarrow T'_3 \rightarrow T'_1$ 
induced by $\tau$-conjugates of $\phi_{1,2} , \phi_{2,3} ,\phi_{3,1}$. 

Let $G= Aut(A)$. 
Let $G^+$ be the subgroup of $G$ generated by all edge-stabilisers. 
It is proved in \cite{tits2} that every element of $G^+$ is either 
a rotation or a translation by an even distance. 
By Theorem 5 of \cite{tits2} $G^+$ is simple of finite width.  

Let a rotation $g_0$ stabilise $\{ a_0 ,a_1\} \cup T_4 \cup...\cup T_{k-1}$ 
pointwise and act on the remaining part of $A\setminus T_0$ as follows. 
\begin{quote} 
$g_0 (x )=\phi_{i,i+1} (x)$ for $x\in T_i$ , $i<3$; \\ 
$g_0 (x )=\phi_{3,1}(x)$ for $x\in T_3$.  
\end{quote} 
We assume that $g_0$ act on $T_0$ in the same way for 
the decomposition $T_0 = \{ a_0 \} \cup T'_1 \cup ...\cup T'_{k-1}$, 
i.e. $g_0$ coincides with $\alpha$ on $T_0$. 
In particular  $|g_0 | = 3$. 
As we already know $G^+ =\langle X \rangle$, where $X$ is the conjugacy class of $g_0$. 
Note that $X^{-1} =X$. 

\begin{prop} 
Under the assumptions above the conjugacy class $X= g^{G^+}_0$ 
is path connected in $Aut(A)^+$.  
\end{prop}

{\em Proof.} 
Let $g_1 \in G^+$ belong to the conjugacy class of $g_0$: 
$g_1 =h g_0 h^{-1}$.  
Let $e_1 = h(e_0 )$, i.e. $e_1$ is fixed by $g_1$ and moreover 
$g_1$ fixes $k-3$ neighbours of each endpoint of $e_1$ 
and acts on the remaining three edges by a cycle.  
Find a path between $e_0$ and $e_1$.  

We start with the case when $e_0 =(a_0 ,a_1 )$, 
$e_1 =(a_1 ,b_1 )$ and    
$g_1$ is defined as follows. 
Let 
$$ 
g_1 (b_1 )= b_1 \mbox{ , } g_1 (b_2) =a_0 \mbox{ , } g_1 (a_0 )=b_3 \mbox{ , } g_1 (b_3 )= b_2 ,  
$$  
and $g_1$ take $T_2$ to $T_0$ by $\phi^{-1}_2$, 
$T_0$ to $T_3$ by $\phi_3 \alpha$ and 
$T_3$ to $T_2$ by $\phi_2 \alpha^{-1} \phi^{-1}_2 \phi_{3,2}$.    
We also assume that $g_1$ acts on $T_1$ by $\phi_1 \alpha \phi^{-1}_1$. 
It is easy to see that $g_1$ is of the same conjugacy class with $g_0$. 

Since $g_0$ acts on $T_0$ by $\alpha$, straightforward computations show that 
$g_0 \cdot g_1$ belongs to the conjugacy class of $g_0$ 
where the edge $(a_1 ,b_3 )$ corresponds to $e_0$. 
\begin{quote} 
Let us fix a pair $(h_0 ,h_1 )$  the conjugacy class of the pair $(g_0 ,g_1 )$ as above. 
\end{quote} 
{\em Remark.} 
It is worth noting here that given $g_0$ and $e_1$ 
the element $g_1$ is not a unique element which fixes $e_1$ 
and together with $g_0$ represents the conjugacy class of 
$(h_0 ,h_1 )$.  
In fact for any two isomorphisms $\psi :T_0 \rightarrow T_2$
and $\psi' :T_0 \rightarrow T_3$ satisfying the properties that 
\begin{itemize} 
\item $\psi$- and $\psi'$-conjugates of $\alpha$ on $T_2$ and $T_3$ respectively 
are identified by $g_0$ and 
\item $g^{-1}_0 \psi$- and $g_0 \psi'$-conjugates of $\alpha$ on $T_1$ are the same, 
\end{itemize} 
the automorphism $h$ defined on $T_2$ by $\psi^{-1}$, defined on $T_1$ 
by the $g^{-1}_0 \psi$-conjugate of $\alpha$, defined on $T_0$ by $\psi' \alpha$ 
and on  $T_3$ by $\psi \alpha \psi^{-1} (g_0 )^{-1}$, satisfies the property 
that $(g_0 ,h)$ is of the same conjugacy class with $(h_0 ,h_1 )$. $\Box$

Now assume that for $e_1 =(a_1 ,b_1 )$ the automorphism $g_1$ is defined 
so that $(g_0 ,g_1 )$ is not of the conjugacy class of $(h_0 ,h_1 )$. 
In this case we find an edge $e_2 \not= e_1$ from the point $b_1$ 
so that for some $g_2$ fixing $e_2$, the pair $(g_1 ,g_2 )$ 
belongs to the conjugacy class of $(h_0 ,h_1 )$. 
We now transfom $g_1$ into some $g'_1$ fixing $e_1$ so that 
the pairs $(g_0 ,g'_1 )$ and $(g'_1 ,g_2 )$ 
are of the conjugacy class $(h_0 ,h_1 )$. 
This will give an $X$-path: $g_0$, $g'_1$, $g_2$, $g_1$ 
(because $g^{-1}_2 g^{-1}_1$ will be the converse of $g_1 g_2$, 
i.e. of the same conjugacy class with $g_0$). 

Let $e_2 =(b_1 ,c )$ and $T_c$ be the subtree 
with the root $c$ without $e_2$. 
We will define the action of $g'_1$ so that 
$a_0 \rightarrow b_3 \rightarrow b_2 \rightarrow a_0$ and 
$g_2 (a_1 )\rightarrow c \rightarrow g^{-1}_2 (a_1 ) \rightarrow g_2 (a_1 )$
are $g'_1$-cycles (in particular to guarantee that the sets 
$g^{-1}_0 (T_1 )$, $T_0$ and $g_0 (T_1 )$ are in the same $g'_1$-cycle).  
The definition is by induction on $n$-balls of $a_1$ 
with $n=3,5,7,...$. 
In fact we want to satisfy the condition of remark above 
for $\alpha$ and the action of $g_2$ on $T_c$ simulteneously 
(for appropriate $\psi$ and $\psi'$). 

At the first step define $g'_1$ on the neighbours of $g_2 (a_1 )$ 
so that the cycle $g_2 (a_0 ) \rightarrow g_2 (b_3 ) \rightarrow g_2 (b_2)$ 
is taken to the $g_2$-cycle of neighbours of $c$ (we have three possibilities). 
This together with the assumptions that 
$g_2 (a_1 )\rightarrow c \rightarrow g^{-1}_2 (a_1 ) \rightarrow g_2 (a_1 )$
is a $g'_1$-cycle 
and $(g'_1 )^{-1}$ takes the triple $g_2 (b_3 ), g_2 (b_2) , g_2 (a_0 )$ to 
$g^{-1}_2 (a_0 ) , g^{-1}_2 (b_3 ), g^{-1}_2 (b_2)$
already defines $g'_1$ on neighbours of $c$ because they must go to 
$\{ g^{-1}_2 (a_0 ), g^{-1}_2 (b_3 )$, $g^{-1}_2 (b_2) ,...\}$
so that their $g'_1$-cycles are 3-element.  

Since we want the pair $(g_0 ,g'_1 )$ to be of the type $(h_0 ,h_1 )$ 
now $g'_1$ must take the triple 
$g^{-1}_0 (g_2 (a_1 )), g^{-1}_0 (c), g^{-1}_0 (g^{-1}_2 (a_1 ))$ 
to the triple $g_0 (g^{-1}_2 (a_1 )), g_0 (g_2 (a_1 )), g_0 (c)$ 
(with the natural extension to the next level).  
Then we define $(g'_1)^{-1}$ and $g'_1$ on the subtree of $T_0$ of depth 2 
satisfying the requirement that the $g_0$-conjugate and the $g^{-1}_0$-conjugate 
(defined on $T_2 = g_0 (T_1 )$ and $T_3 = g^{-1}_0 (T_1)$ respectively) 
of the action of $g'_1$ on the subtree of $T_1$ of depth 2 must 
correspond by $g'_1$ to the action of $g_0$ on $T_0$. 
For example we have only three possibilities for the definition 
of $g'_1$ on the non-trivial $g_0$-cycle in $T_0$.  
This finishes the definition of $g'_1$ on the 3-ball of $a_1$. 

The next step (for the 5-ball) is similar. 
Taking the $g_2$-conjugate (on $g_2 (A\setminus T_1 )$) 
of the action of $g'_1$ on the subtree of depth 3 in $A\setminus T_1$ 
we extend $g'_1$ to the next two levels of $g_2 (A\setminus T_1 )$ 
so that it corresponds to the $g_2$-action on $T_c$. 
Then as in the previous paragraph we extend $g'_1$ to the 
subtrees of $T_1$, $T_3$ and $T_0$ of depth 4. 
This defines $g'_1$ on the 5-ball of $a_1$. 
 Continuing this procedure we define the required $g'_1$. 
This finishes the case when $e_0$ and $e_1$ have one common endpoint. 
 
The general case is obviously reduced to the situation when  
$e_0$ and $e_1$ do not have common endpoints. 
Then we can build an $X$-path of between $g_0$ and $g_1$ as follows. 
Take the shortest path $e_0 ,e'_1 ,...,e'_k, e_1$ connecting $e_0$ and $e_1$ 
and applying the argument above build an appropriate sequence 
$g_0 ,g'_1 ,...., g'_l ,g_1$ (with $k\le l$) of conjugates of $g_0$ 
so that at every step they realise a conjugate of $(h_0 ,h_1)$.   
$\Box$

\section{Median pretrees and actions of groups}

In this section we develop some basic tools which we will
apply to the main theorem of the paper (Theorem \ref{main} ).
Here we consider a very general situation when a group acts on a median pretree. 

\subsection{Median pretrees} 

The following definition is taken from \cite{bow}, p.10. 
It is related to the definition of a B-relation given in \cite{an}, p.55.

\begin{Def}
A ternary structure $(T,B)$ is a pretree if the following axioms
are satisfied: \parskip0pt

* $(\forall x,y)(\neg B(y;x,x)$

* $(\forall x,y,z)(\neg(B(y;x,z) \wedge B(z;x,y)))$;

* $(\forall x,y,z)(B(y;x,z) \leftrightarrow B(y;z,x))$;

* $(\forall x,y,z,w)(B(z;x,y) \wedge z \not= w \rightarrow
(B(z;x,w) \vee B(z;y,w)))$. 
\end{Def}

Define open and closed {\em intervals with endpoints} $t,t'$ as follows: 
$$
(t,t') = \{ x \in T: B(x;t,t') \} \mbox{ and } 
[t,t'] = \{ x \in T: B(x;t,t') \mbox{ or } x = t \mbox{ or } x = t'\} . 
$$
The intervals $[t,t')$ and $(t,t']$ are defined similarly. 
We also call endpoints {\em extremities} (as in \cite{tits}) 
and closed intervals {\em segments}. 
A nonempty subset $S \subseteq T$ is an {\em arc}, if $S$ is
{\em full} (that is for all $x,y \in S$, $[x,y] \subseteq S$) and
{\em linear} (for all distinct $x,y,z \in S$ we have
$B(y;x,z)$ or $B(z;x,y)$ or $B(x;y,z)$). \parskip0pt

A point $x \in T$ is {\em terminal}, if
$x$ is never between any $y,z \in T$.
The pretree $T$ can be naturally decomposed
$T = T_{0} \cup P$, where $P$ is the set of all terminal 
points and $T_0$ is the set of all non-terminal ones.  \parskip0pt

The pretree $(T,B)$ is {\em median} if for any $x,y,z \in T$
there is an element  $c \in [x,y]\cap [y,z] \cap [z,x]$.
In this case $c$ is unique and is called the {\em median}
of $x,y,z$; we will write $c = m(x,y,z)$. 
In Remark 2 of \cite{ivalox} some connections 
with median metrics are described. \parskip0pt

It is clear that every simplicial or real tree (see Definition \ref{RT}) 
can be considered as a {\em complete} median pretree 
(i.e. every arc is an interval, not necessarily closed). 
We may also add ends as terminal points.

\begin{Assumption}
{\em From now on we consider only median pretrees.}
\end{Assumption} 

Let $x \in T_{0}$.
A maximal arc of the form $L_{x} =\bigcup [x,t_{\gamma}]$
where all $t_{\gamma}$ are not terminal, is called 
a {\em half-line} starting at $x$.
Half-lines $L_{x}$ and $L_{y}$ are equivalent
if there is a half-line $L_{z} \subseteq L_{x} \cap L_{y}$.
An {\em end} ${\bf e}$ of $T$ is an equivalence class of
half-lines.
Define a partial order $<_{\bf e}$ by $x <_{\bf e} y$
if the half-line $L_{y}$ representing ${\bf e}$ contains $x$.
Sometimes we will write $L_y =(-\infty ,y]$. 
\parskip0pt

It is clear that an arc of the form $[x,p)$, where $p$ 
is terminal ($\in P$), is a half-line.
By the definition terminal points, any pair of half-lines $[x,p)$
and $[y,p)$ with $x,y\in T_0$, have a common point from $T_{0}$
(which is the corresponding median).
This shows that the set of all half-lines $[t,p), t\in T_{0}$
forms an end (the {\em end corresponding} to $p\in P$).\parskip0pt

A maximal arc of the form $\bigcup [t_{\gamma},t'_{\gamma}]$
where $t_{\gamma},t'_{\gamma}$ are not terminal, is called a
{\em line}.
It is worth noting that the ends of a line of a complete pretree $T$ 
are presented by a pair of termial points of $T$.

The following lemma is conceivably known (see \cite{tits},\cite{bc}
and Section 2 of \cite{bow}) and is based on existence of medians.
A complete proof of the lemma (in a slightly more general situation)
is given in \cite{ivanov}.

\begin{lem} \label{4sts}
(1) The intersection of a finitely many (distinct) segments or
half-lines with a common extremity $t \in T_{0}$ (= non-terminal points) 
is a segment having $t$ as an extremity. \parskip0pt

(2) If $t,q,r \in T_{0}$ satisfy $[t,q] \cap [q,r] = \{ q\}$, then
$[t,r] = [t,q] \cup [q,r]$.

(3) If the interval $[t,q)$ is not closed and $[t,q)=[t,r)$, then
$q = r$. \parskip0pt

(4) Let each of $A$ and $B$ be a segment or a line and $|A \cap B|\le 1$.  
Then there exists $t\in A$  and $q\in B$ so that 
the segment $[t,q]$ is contained in every full set 
which has a non-empty intersection with both $A$ and $B$. 
The points $t$ and $q$ are determined uniquely. 

(5) Let $T$ be complete and $A$ and $B$ be two full subsets of $T_{0}$ 
whose intersection consists of at most one point.
Then there exists a segment $[t,q]$ which is contained in every
full set which has a non-empty intersection with both $A$ and $B$,
and, moreover, $A\cup \{ t\}$ and $B\cup \{ q \}$ are full.
If for some $\varepsilon,\tau \in \{ 0,1\}$, $t$ and $q$ satisfy
$t\in^{\varepsilon} A \wedge q\in^{\tau} B$ (where $\in^{0}$ denotes
$\not\in$), then the condition $t\in^{\epsilon} A \wedge q\in^{\tau} B$ 
determines the segment uniquely. 
\end{lem}

We just mention how the point $t$ in 
Statement 4 can be found. 
If $|A \cap B|=1$, then $\{ t\} = A\cap B$. 
Let $A\cap B=\emptyset$. 
If $A=[w,v]$, take any $b\in B$ and let $t:=m(b,w,v)$.   
If $A$ is a line take any $b\in B$ and find $w,v\in A$ 
so that $\{ b,w,v\}$ is not linear. 
The existence of such a pair follows from 
the maximality condition in the definition of a line. 
Now take $t:= m(b,w,v)$. 
The point $q$ can be found similarly. 

Let $T$ and $A,B\subseteq T_{0}$ with $A\cap B=\emptyset$, satisfy 
the assumptions of Lemma \ref{4sts} (4 or 5). 
Then we call the interval $[t,q]\setminus (A\cup B)$ the
{\em bridge between} $A$ and $B$.
It can happen that the bridge is an open (or empty) interval.
In this case we define the bridge by its extremities $t$ and $q$ 
as $(t,q)$.

\subsection{Flows}

We now describe one of our tools.
A binary relation $r$ (a partial ordering, where
$r(a,b)\vee (a=b)$ is interpreted as $a\le b$) on a pretree 
$T$ is called a {\em flow} (\cite{bow}, pp. 23 - 25) if
it satisfies the following axioms: \parskip0pt

* $\neg (r(x,y)\wedge r(y,x))$,

* $B(z;x,y)\rightarrow r(x,z)\vee r(y,z)$, 

* $(r(x,y)\wedge z \not=y) \rightarrow (B(y;x,z)\vee r(z,y))$.

The material of the next paragraph is based on pp. 26 - 28 of 
\cite{bow}. 
It would be helpful for the reader (but not necessary) 
to recall some formulations given there. \parskip0pt

By Lemma 3.8 of \cite{bow} for any endless directed arc $(C,<)$ 
the formula 
$$
(x,y)\in r\leftrightarrow\exists z\in C\forall w >z B(y;x,w)  
$$
defines a flow. 
In this case we say that the flow $r$ is
{\em induced by an endless directed arc} $(C,<)$ (see \cite{bow}, p. 26).  
Then it is easy to see that for any arc $J$ if $J$ does not have
maximal elements with respect to $r$, then the formula
$$ 
E^{J}_{r}(x,y) = r(x,y)\vee r(y,x)\vee x = y
$$ 
defines an equivalence relation on $J$ with at most two classes 
(apply Lemma \ref{4sts}(4) to two lines containing $C$ and $J$ respectively).
\begin{Def} 
We say that the flow $r$ {\em lies} on $J$ if 
\begin{itemize} 
\item $J$ contains a maximal element of $r$ in $T$ or 
\item $J$ does not have any element which is the greatest in $J$ 
with respect to $r$ and the equivalence relation
$E^{J}_{r}(x,y)$ defines a non-trivial cut
$J= J^{-}\cup J^{+}$, such that for any $a \in J^{+}$,
$b\in J^{-}$ there is no $c \in T$ with $r(a,c)\wedge r(b,c)$
(so it may happen that $C\cap J$ is cofinal with $C$).
\end{itemize} 
\end{Def} 
It is easy to see that for any line $L$ there is a natural 
function from the set of all flows of $T$ induced by endless 
directed arcs and lying on $L$ onto the set of all Dedekind 
cuts on $L$: the Dedekind cut corresponding to a flow $r$ is 
determined by a maximal element or by the equivalence relation 
$E^{J}_{r}(x,y)$.
The following lemma shows that when $T$ is median and dense, 
this correspondence is bijective for flows without maximal elements.

\begin{lem} \label{flowcut} 
Assume that $T$ is a median pretree and $r$ is a flow induced by 
an endless directed arc. 
Let $C,D$ be arcs of $T_0$ which are linearly ordered with respect 
to $r$ and do not have upper $r$-bounds in $T_0$.
Then $C$ and $D$ are cofinal or for any $c \in C$ and $d \in D$,
$[c,d]=\{ a: a\in C\wedge r(c,a)$ or $a\in D\wedge r(d,a)\}$.
Each of these arcs defines the flow $r$ as in the definition above.
\end{lem}

{\em Proof.}
To see the first statement we start with the case when 
$r$ is induced by $C$.
Let $c$ and $d$ be as in the formulation. 
If $C$ and $D$ are not cofinal then $C\cap D=\emptyset$
(if $t\in C\cap D$ and $D\models t<t'$,
then $t'$ must belong to $C$).
By the definition of $r$, if $D\models d<d'$, then $d'$
belongs to some $[d,c']$ with $C\models c<c'$.
If $d'\not\in [d,c]$, then $d'$ belongs to $[c^{*},c']$, where
$c^{*}$ is the median of $c,c',d$, a contradiction with
$C\cap D=\emptyset$.
It is now easy to see that the set $[c,d]\cap D$ consists 
of all $d'\in D$ with $r(d,d')$. \parskip0pt

This implies that there is no $a\in [c,d]\setminus (C\cup D)$
(otherwise $a$ is an upper $r$-bound for $D$) and
there is no $c'\in C$ with $c'\not\in [c,d]$ and $c<c'$
(otherwise $c'$ is an upper $r$-bound for $D$).
Now a straightforward argument gives the formula for $[c,d]$ 
as above. \parskip0pt

In the case when $r$ is induced by some ordering $A$ and
$A$ is cofinal with $C$ or $D$, the argument above works again.
If $A$ is not cofinal with these orderings then for every $a\in A$
we have
$[c,a]=\{ a': a'\in C\wedge r(c,a')$ or $a'\in A\wedge r(a,a')\}$
and
$[a,d]=\{ a': a'\in D\wedge r(d,a')$ or $a'\in A\wedge r(a,a')\}$.
Then the median of $a,c,d$ must belong to one of the intervals $C$, 
$D$ or $A$. 
The case $m(a,c,d)\in A$ is impossible, because the elements of $A$ 
greater than $m(a,c,d)$ cannot belong both to $[c,a]$ and $[d,a]$. 
When $m(a,c,d)\in C\cup D$, the arcs $C$ and $D$ are cofinal. 
\parskip0pt

To show that {\em the flow $r$ is induced by
any of its linear orderings without upper $r$-bounds in $T_0$}
we apply similar arguments as above.
Indeed, if $A$ and $C$ are linear orderings without upper $r$-bounds
in $T_0$, $r$ is induced by $A$ and $C$ is not cofinal with $A$,
then since $T$ is median any $D$ as above is cofinal with $A$ or $C$.
If $D$ is cofinal with $A$, then $D$ induces $r$.
If $D$ is cofinal with $C$, then for any $d\in D$ and $a\in A$,
$[a,d]=\{ a':a'\in D\wedge r(d,a')$ or $a'\in A\wedge r(a,a')\}$.
Let $t\in T_0$ and $t^{*}=m(t,a,d)$ for some $d\in D$ and $a\in A$.
Now it is straightforward that for any $t'\in T_0$ the condition
$r(t,t')$ is equivalent to the disjunction:  
$t'\in (t,t^{*}]$ or $t'\in (t^{*},d]\cap A$ or $t'\in (t^{*},a]\cap D$.
Using this formula it is easy to verify that $r$ is induced by $D$ 
(notice that $A$ and $D$ are symmetric in this condition). $\Box$ 
\bigskip 

By Lemma \ref{flowcut} we see that when $r$ is a flow of 
a dense median pretree $T$ defined by an endless directed arc  
and $r$ lies on a line $L$ but does not have a maximal element 
lying on $L$, then each of the half-lines on $L$ defined by the 
corresponding equivalence relation is an endless directed arc 
inducing $r$.

\subsection{Actions of groups} 

Let $G$ be a group acting on a median pretree
$T$ by automorphisms of the structure $(T,B)$.
It is clear that the set $P$ of terminal points is $G$-invariant.
Let $g \in G$.
The set of $g$-fix points is denoted by $T^{g}$.
We say that $g$ is elliptic if $T^g_0 \not=\emptyset$, 
i.e. $g$ does not have non-terminal fixed points. 
The element $g\in G$ is {\em loxodromic}, if 
$T^{g}_{0} =\emptyset$, $|g| = \infty$ and there exists a unique
$g$-invariant line in $T_{0}$ such that $g$ preserves
the natural orders on the line.
It is called the {\em axis (characteristic line)} of $g$.
In the case of an isometric action on an $\mathbb{R}$-tree
a loxodromic element is {\em hyperbolic}.
If $g$ is elliptic or loxodromic, then we denote by $C_g$ 
the set $T^g$ or $L_g$ respectively. 
We will say that $C_g$ is the {\em characteristic set} of $g$. 

The following lemma is Lemma 1.5 in \cite{ivanov}. 
The proof is based on standard arguments from \cite{tits} (Section 3.1).

\begin{lem} \label{loxo}
Let $g$ be an automorphism of a median pretree $T$ and 
$T_0$ be the set of all non-terminal points of $T$.
If $g$ is loxodromic, then 
\begin{enumerate} 
\item  for any $p \in T_{0}$ the segment
$[p,g(p)]$ meets the characteristic line $L_{g}$
and $[p,g(p)] \cap L_{g} = [q,g(q)]$ for some $q \in L_{g}$.
\item $x \in L_{g}$ if and only if $x$ is the median of
$x,g^{-1}(x),g(x)$.
\end{enumerate} 
\end{lem}

The following proposition is a very particular case 
of Theorem 6 in \cite{ivalox}. 
It can be considered as a metric-free 
version of statements 1.5 and 1.6 of \cite{cullmorg}.
Since the proof is easy we include it for completeness. 

\begin{prop} \label{3loxo}
Let $G$ act on a median pretree $T$.
Let $h_{1},h_{2},h_{3} \in G$ be loxodromic and
$h_{2} \cdot h_{1} = h_{3}$.
Then 
if $\{ i,j,k\} =\{ 1,2,3\}$ and $L_{h_i} \cap L_{h_j}=\emptyset$ 
then the bridge from $L_{h_i}$ to $L_{h_j}$ belongs to 
$L_{h_k}$. 
\end{prop} 

{\em Proof.} 
Since 
$h_3 = h_2 \cdot h_1 \Leftrightarrow h_1 = h^{-1}_2 \cdot h_3 \Leftrightarrow h_2 = h_3 \cdot h^{-1}_1$ 
and the axis of a loxodromic element $h$ coincides with the axis 
of $h^{-1}$, it suffices to consider the case $L_{h_1} \cap L_{h_2}= \emptyset$. 
If $a\in L_{h_1}$, $b\in L_{h_2}$ and $(a,b)$ is the bridge between $L_{h_1}$ 
and $L_{h_2}$ (see Lemma \ref{4sts}(4)), then by Lemma \ref{loxo}(1) 
$L_{h_3}$ intersects both $[h^{-1}_1 h^{-1}_2 (a), a]$ and $[b, h_2 h_1 (b)]$. 
Thus it contains $[a,b]$. 
$\Box$

\section{Non-nesting actions on real trees} 

In Section 2.1 we collect well-known facts about non-nesting actions 
on $\mathbb{R}$-trees (\cite{bc}, \cite{guiva} and \cite{levitt}). 
In Section 2.2 we prove some basic lemmas 
which will be applied in the proof of Theorem \ref{main}.

\subsection{Non-nesting actions on $\mathbb{R}$-trees}

\begin{Def} \label{RT} 
An $\mathbb{R}$-tree is a metric space $T$ such that for any $x\neq y\in T$,
there is a unique topologically embedded arc joining $x$ to $y$, and this arc 
is isometric to some interval of $\mathbb{R}$.
\end{Def} 

We define $[x,y]$ as the arc joining $x$ to $y$ 
if $x\neq y$, and $[x,y]=\{x\}$ if $x=y$. 
We say that $[x,y]$ is a \emph{segment}.  
A subset $S\subseteq T$ is {\em convex} if 
for any $x,y \in S$, $[x,y] \subseteq S$. 
A convex subset is also called a {\em subtree}. 

Define a betweenness relation by $B(x;y,z) \Leftrightarrow x\in (y,z)$. 
Now $T$ can be considered as a pretree. 
It is easy to see that $(T,B)$ is a median pretree without terminal points. 

Any homeomorphism of $T$ clearly preserves the betweenness relation. 
Sometimes automorphisms of $(T,B)$ are called  
{\em weak homeomorphisms} of the $\mathbb{R}$-tree $T$. 
All actions on $T$ are via weak homeomorphisms.

\begin{Def}
An action of $G$ on a pretree $T$ by its automorphisms is {\em non-nesting} if there 
is no segment $I\subseteq T$, and no $g\in G$ such that $g(I)\subsetneqq I$. 
\end{Def}

From now on, we assume that $G$ has a non-nesting action on 
an $\mathbb{R}$-tree $T$ by weak homeomorphisms.
The following lemma in particular says that when $g\in G$ 
is not {\em elliptic}, it is {\em loxodromic}.

\begin{lem}[{\cite[Theorem 3]{levitt}}]\label{lem_nnest}
Let $G$ be a group with a non-nesting action on an $\mathbb{R}$-tree $T$.

  \begin{itemize}
  \item If $g$ is elliptic, its set of fix points $T^g$ is a closed convex (= full ) subset.
  \item\label{nnest_it2} If $g$ is not elliptic, there exists a unique line $L_g$ preserved by $g$; 
  moreover, $g$ acts on $L_{g}$ by  an order preserving transformation, 
  which is a translation by a real number up to topological conjugacy. 
  \end{itemize} 
\end{lem}

In \cite{levitt}, $g$ is assumed to be a homeomorphism, 
but the argument still applies, see Lemma 1.5 in \cite{guiva}. 
When $g$ is loxodromic, the action of $g$ on $L_g$ defines 
a natural ordering on $L_g$ such that for all $x\in L_g$, $x<g(x)$. 
Then by Lemma \ref{lem_nnest} the element $g$ is strictly increasing on $L_{g}$.

This can be complemented as follows. 
Let $L$ be a line of $T$ and 
$G_{\{ L\} }$ (and $G_{L}$) be the stabilizer
(pointwise stabilizer) of $L$ in $G$.
Assume that $G_{\{ L\} }$ does not have elements reversing the ends of $L$
(this happens when it does not have subgroups of index 2).
Then non-nesting implies that given an ordering of $L$ 
corresponding to the betweenness relation of $T$, the following
rule makes $G_{\{ L\} }/G_{L}$ a linearly ordered group:
$g<g'$ if and only if there exists $t\in L$ so that $g(t)< g'(t)$ 
(equivalently for all $t\in L$, $g(t) <g'(t)$). \parskip0pt

Since up to topological conjugacy a loxodromic $g$ 
can be viewed as a translation by a real number, we see that 
{\em for every $h\in G_{\{ L_{g}\} }$ there exists $n\in \omega$ such that}
$h<g^{n}$. 

\begin{Remark} 
If $H$ is a subgroup of the stabiliser$G_{\{ L\} }$ then 
by completeness $L$ consists of $H$-invariant intervals.
If $[a,b)$ is such an interval then by non-nesting,
$h(a) \not\in [a,b)$, where $h\in H$.
We see that $H$ acts trivially on $L$ or
there are no proper $H$-invariant arcs in $L$.
\end{Remark} 

The following statement is Proposition 1.7(a,b) in \cite{guiva}. 

\begin{prop} \label{2.6}
Let $G$ be a group with a non-nesting action on an $\mathbb{R}$-tree $T$.
\begin{enumerate}
\item \label{it1} If $g$ is elliptic and $x\notin T^g$, then $[x,g(x)]\cap T^g=\{a\}$ where 
$a$ is the projection of $x$ on $T^g$. 
\item\label{it2}
If $g,h\in G$ are elliptic and $T^{g}\cap T^h=\emptyset$,
then $gh$ is loxodromic, its axis contains the bridge between $T^g$ and $T^h$,
and $T^g\cap L_{gh}$ (resp.\ $T^h\cap L_{gh}$)
contains exactly one point.  
In particular, if $g,h$ and $gh$ are elliptic, then $T^g\cap T^h\cap T^{gh}\neq\emptyset$. 
\end{enumerate}
\end{prop}

\subsection{Subgroups stabilising an end} \label{*action}

Let a group $G$ have a non-nesting action on
an $\mathbb{R}$-tree $T$ without $G$-fixed points.
Assume that there is a loxodromic $g \in G$. 
Let $a_{0} \in L_{g}$ and ${\bf e}$ be the end represented by
$(-\infty ,a_{0}]$ ($-\infty$ is chosen so that $g$ is increasing).
Consider the stabiliser $G_{{\bf e}}$.
Let $G_{({\bf e})}$ be the subset of $G_{{\bf e}}$ of all
elements fixing some points. 
Note that any $h \in G_{{\bf e}}$ defines a map
$(-\infty ,a] \rightarrow (-\infty ,b]$ for some $a,b \le a_{0}$.
Thus $G_{({\bf e})}$ is the subgroup of elements fixing pointwise
half-lines $(-\infty ,a]$ for some $a \le a_{0}$ (by non-nesting).
We also see that it is normal in $G_{{\bf e}}$.
In the following lemma we consider the group
$G_{{\bf e}}/G_{({\bf e})}$.

\begin{lem} \label{1.10}
Let a group $G$ have a non-nesting action on an
$\mathbb{R}$-tree $T$ without $G$-fixed points.
Let $g \in G$ be loxodromic and ${\bf e}$ be the
$(-\infty)$-end of $L_{g}$.
Then the group $G_{{\bf e}}/G_{({\bf e})}$ is embeddable into
$(\mathbb{R},+)$ as a linearly ordered group under the ordering: 
$gG_{({\bf e})}\prec g'G_{({\bf e})}\Leftrightarrow \exists t\in T (g'(t)<_{{\bf e}}g(t))$.
\end{lem}

{\em Proof.} 
By Lemma \ref{lem_nnest} any $h \in G_{{\bf e}}\setminus G_{({\bf e})}$ 
defines a map $(-\infty ,a] \rightarrow (-\infty ,b]$ with some 
$a,b\in (-\infty ,a_0 ]$, which is strictly monotonic on $(-\infty ,a]$.
Now it is easy to see that the ordering $\prec$ of the set 
$G_{{\bf e}}/G_{({\bf e})}$ is a linear. 
To see that $G_{{\bf e}}/G_{({\bf e})}$ is a linearly ordered group 
note that if $g_{1},g_{2} \in G_{{\bf e}}$ satisfy $g_{1}(a) < g_{2}(a)$
for some $a \le a_{0}$, then by non-nesting
for all $a' \le a$, $g_{1}(a') < g_{2}(a')$.
Thus for every $g' \in G_{{\bf e}}$ there exists $b \le a$
such that $g_{1}\cdot g'(x) \le g_{2}\cdot g'(x)$ for all
$x \le b$.
On the other hand if $g_{1},g_{2} \in G_{{\bf e}}$ satisfy
$g_{1}G_{({\bf e})} \preceq g_{2}G_{({\bf e})}$ then obviously
$g'\cdot g_1 G_{({\bf e})}\preceq g'\cdot g_{2} G_{({\bf e})}$ 
for all $g' \in G_{{\bf e}}$.
This shows that $G_{{\bf e}}/G_{({\bf e})}$ is a linearly
ordered group. \parskip0pt

Since the elements of $G_{{\bf e}}$ act by translations up to
topological conjugacy, $G_{{\bf e}}/G_{({\bf e})}$ is Archimedean.
By H\"{o}lder's theorem it is a subgroup of $(\mathbb{R},+)$
\cite{birkhoff}. 
$\Box$

\bigskip

Let $g$, ${\bf e}$ and $a_0$ be as above. 
{\em We now define an action $*_{g}$ of} $G_{{\bf e}}$
on $L_{g}$.
Let $h\in G_{{\bf e}}$ and $c \in L_{g}$.
Find a natural number $n_{0}$ such that $g^{-n_{0}}(c)$ is 
less with respect to $<_{{\bf e}}$ than any of 
the elements $a_0 ,h(a_0 ), h^{-1}(a_0 )$.
Now let $h*_{g}c = g^{n_{0}}hg^{-n_{0}}(c)$.
By the choice of $n_{0}$ we see 
$h*_g c \in L_{g}$. \parskip0pt

It is worth noting that for every $n \ge n_{0}$,
$h*_{g}c = g^{n}hg^{-n}(c)$.
This follows from the fact
that the element $h^{-1}g^{n-n_{0}}hg^{n_{0}-n}$ belongs to
$G_{({\bf e})}$ (as $G_{{\bf e}}/G_{({\bf e})}$ is a subgroup of
$(\mathbb{R},+)$) and then (by non-nesting) the transformations
$hg^{n_{0}-n}$ and $g^{n_{0}-n}h$ are equal at $g^{-n_{0}}(c)$.
We now see:
$$
g^{n}hg^{-n}(c) =
g^{n_{0}}h(h^{-1}g^{n-n_{0}}\cdot hg^{n_{0}-n}(g^{-n_{0}}(c)))
=g^{n_{0}}hg^{-n_{0}}(c).
$$
Now it is easy to see that $*_{g}$ is an action and the elements
of $G_{({\bf e})}$ act on $L_{g}$ trivially. \parskip0pt

Let $L^{a_{0}}_{g}= G_{{\bf e}}a_{0} \cap L_{g}$ where $G_{{\bf e}}a_{0}$ 
is the orbit of $a_{0}$ under the action of $G_{{\bf e}}$ on $T$.
Then there exists a surjection 
$\nu_{a_{0}} : G_{{\bf e}}\rightarrow L^{a_{0}}_{g}$ defined
by $\nu_{a_{0}} (h) = h*_{g}a_{0}$ (with respect to the
action defined above).
It is easy to see that the map $\nu_{a_{0}}$ is surjective.
Moreover, for any $h,h'\in G_{{\bf e}}$,
$\nu_{a_{0}} (h\cdot h') = h*_{g}\nu_{a_{0}} (h')$.

\begin{lem}  \label{1.11}
The map $\nu_{a_{0}}$ defines an order-preserving bijection
from $G_{{\bf e}}/G_{({\bf e})}$ onto $L^{a_{0}}_{g}$ under
the order induced by $L_{g}$.
\end{lem}

{\em Proof.}
Notice that if $\nu_{a_{0}} (h_{1}) = \nu_{a_{0}} (h_{2})$
then $h_{1}h^{-1}_{2}$ fixes some $(-\infty ,a]$,
$a \le a_{0}$, pointwise.
Indeed, let $n$ and $a = g^{-n}(a_{0})$ be chosen so
that $a,h_{1}(a),h_{2}(a)$, $h^{-1}_{1}(a),h^{-1}_{2}(a)$
$\in L_{g} \cap L_{h_{1}} \cap L_{h_{2}}$
(for $h_i \in G_{({\bf e})}$ we replace $L_{h_i}$ by $T^{h_i}$).
Then
$$
g^{n}h_{1}g^{-n}(a_{0}) = \nu_{a_{0}} (h_{1}) =
\nu_{a_{0}} (h_{2}) = g^{n}h_{2}g^{-n}(a_{0})
$$
and we see that $h_{1}(a) = h_{2}(a)$.
Now the claim follows from non-nesting. \parskip0pt

Applying non-nesting again we obtain that the preimage of 
$a_{0}$ (with respect to $\nu_{a_{0}}$) equals the subgroup
$G_{({\bf e})}$ of elements fixing pointwise $(-\infty ,a]$ 
for some $a \le a_{0}$. \parskip0pt

The proof of Lemma \ref{1.10} shows that the condition
$h_{1} G_{({\bf e})}\prec h_{2} G_{({\bf e})}$ means the existence 
of $a \in L_{g}$ with $h_{1}(a') < h_{2}(a')$ for all $a'<a$.
This obviously implies $\nu_{a_{0}}(h_{1})< \nu_{a_{0}}(h_{2})$.
We see that $(L^{a_{0}}_{g},<)$ can be identified with
the group $(G_{{\bf e}}/G_{({\bf e})},\prec )$. $\Box$

\begin{lem} \label{1.12}
If the ordering of $L^{a_{0}}_{g}$ is not dense, then
$L^{a_{0}}_{g}$ is a cyclic group with respect to
the structure of $G_{{\bf e}}/G_{({\bf e})}$.
\end{lem}

{\em Proof.} Notice that if there is an interval $(a,b)$,
$a,b \in L^{a_{0}}_{g}$, which does not have elements from
$L^{a_{0}}_{g}$, then every $c\in L^{a_{0}}_{g}$ has
a successor from $L^{a_{0}}_{g}$.
This follows from the fact that $[a,g(a)]$ and $[c,g(c)]$ can
be taken onto $[g^{-1}(a_{0}),a_{0}]$ by an element from
$G_{{\bf e}}$.
$\Box$

\section{Equivariant embeddings}

\subsection{A theorem on equivariant embedding} 

Consider the situation of the previous section. 
Let a group $G$ have a non-nesting action on 
an $\mathbb{R}$-tree $T$ without $G$-fixed points in $T$.
Let ${\bf e}$ be an end and $(-\infty ,a_{0}]$ be a half-line 
representing it. 
Assume that there is a loxodromic element in the stabiliser $G_{{\bf e}}$ 
so that $(-\infty ,a_0 ]$ contains a half-line of its axis. 
Thus the group $G_{{\bf e}}/G_{({\bf e})}$ is non-trivial.

\begin{Def} 
Let $C$ be a subset of $T$. 
We say that $C$ is {\bf e}-{\em contractible} if for any $c\in C$ 
there is $g\in G_{{\bf e}}$ so that $g(c)\in (-\infty ,a_0 ]$. 
\end{Def} 

Note that the definition does not depend on the choice of the 
half-line $(-\infty, a_0 ]$ in ${\bf e}$. 
To illustrate the notion observe that for any loxodromic $g\in G_{{\bf e}}$ 
the line $L_g$ is {\bf e}-contractible. 
When $G$ is the group of all isometries of a regular simplicial tree $T$ 
then considering $T$ as a real tree it is easy to see that 
any subset of $T$ is {\bf e}-contractible with respect to any end {\bf e} 
(any line of $T$ is an axis of some loxodromic element). 
These examples show that the condition of ${\bf e}$-contractibility 
for appropriate ends is very natural for characteristic subsets $C_g$ 
of typical real $G$-trees. 

\begin{Remark} 
Let $g\in G_{{\bf e}}$ be loxodromic. 
Thus $(-\infty, a_0 ]$ has a cofinal intersection with $L_g$.  
Let $T_{{\bf e}} = \bigcup \{ s(L_g ): s\in G_{{\bf e}}\}$. 
Note that $T_{{\bf e}}$ is full and 
if $C$ is ${\bf e}$-contractible, then $C\subseteq T_{{\bf e}}$. 
\end{Remark} 

To formulate our theorem let us introduce the following notion.
For a subset $A$ of a median pretree $T$ define  $c(A)$, 
the {\em median closure} of $A$, as the minimal subpretree 
of $T$ with the property that the $T$-median of any triple from 
$c(A)$ belongs to $c(A)$.
It is clear that in the case when a group $G$ acts on 
a pretree $T$ and $A$ is $G$-invariant, the pretree 
$c(A)$ is $G$-invariant too.  

We also remind the reader that a real $G$-tree $T$ is 
minimal if it does not have proper $G$-invariant subtrees 
(i.e. convex subsets). 
It is clear that any $G$-tree has a minimal subtree.  

In this section we prove the following theorem. 

\begin{thm} \label{main}  
Let a group $G$ have a non-nesting action on an $\mathbb{R}$-tree $T$, 
so that the action is minimal.  
\begin{itemize} 
\item Let $X \subseteq G$ be a path connected conjugacy class with 
$G=\langle X\rangle$,  
\item assume that for any $h\in X$ there is an end ${\bf e}$ so that the 
characteristic set $C_h\subset T$ is {\bf e}-contractible. 
\end{itemize} 
Then for any $v\in T$ the $G$-subpretree $c(Gv) \subseteq T$ 
satisfies one of the following properties: 
\begin{itemize} 
\item for any $h\in X$ and an end ${\bf e}$ contracting $C_h$, any orbit of 
the stabiliser $G_{{\bf e}}$ is not dense in $c(Gv)\cap C_h$  
\item the action of $G$ on $c(Gv)$ can be extended to an isometric action of 
$G$ on a real tree. 
\end{itemize} 
\end{thm}

\subsection{Remarks} 

\begin{Remark} \label{more} 
We will prove a slightly stronger statement. 
Assume that in the situation described by the formulation the point $v$, 
the element $h\in X$ and the end ${\bf e}$ witness that the first possibility of 
the conclusion fails and assume that for some loxodromic $g\in G_{{\bf e}}$ 
and some $a_0 \in L_g$ the $G_{{\bf e}}$-orbit $L^{a_0}_g=G_{{\bf e}}a_{0}\cap L_{g}$ 
is dense (the case when it is discrete is easy).  
Then we prove that {\em any $G_{{\bf e}}$-invariant metric 
(with respect to the action $*_g$) of the Dedekind completion of $L^{a_0}_g$  
can be extended to a metric of some real tree as in the second case of 
the conclusion.}  
\end{Remark} 

\begin{Remark}  
It is worth noting that the condition of ${\bf e}$-contractibility of $C_h$ 
is automatically satisfied when $h$ is loxodromic. 
In this case $C_h$ is an axis and any its end can be chosen as the corresponding ${\bf e}$. 
\end{Remark} 

\begin{Remark} 
To illustrate some aspects of the formulation let $g\in G$ 
be loxodromic and $a_{0}\in L_{g}$.
It is clear that the set $L^{a_{0}}_{g}$ 
is cofinal (in both directions) in the line $L_{g}$.
On the other hand it may happen that the ordering
$c(Ga_{0})\cap L_g$ (induced by a natural ordering of $L_{g}$)
is dense but not dense in $L_{g}$.  

Indeed consider $\mathbb{R}$ as
$\mathbb{Z}\times \{ a,b\}\times (0,1]$, where the elements of
$(2k,2k+1]$ are denoted by triples $(k,a,r)$, $r \in (0,1]$ and
the elements of $(2k+1,2k+2]$ are denoted by triples $(k,b,r)$,
$r \in (0,1]$, $k\in \mathbb{Z}$.
Let $G = \mathbb{Q}$  act on $\mathbb{R}$ as follows.
If $q + k + r = k' + r'$ , with $k'\in \mathbb{Z}$ and
$r' \in (0,1]$, then put $(k,a,r)+q = (k',a,r')$ and
$(k,b,r)+q = (k',b,r')$.
It is easy to see that the action of $\mathbb{Q}$ obtained on
$\mathbb{R}$ is non-nesting.
On the other hand the interval $(0,1]$ (consisting of all
$(0,a,r)$ with $r \in (0,1]$) does not contain any element of
the orbit of $0= (-1,b,1)$.
The fact that the orbit of $0$ is a dense ordering follows from
the density of $\mathbb{Q}$.
Its median closure coincides with the orbit. 
\end{Remark}

\begin{Remark}
We now give an example of a minimal isometric action of a group $G$ 
on a real tree $T$ so that $G$ is generated by a conjugacy class $X$ 
of loxodromic elements and there is a point $v\in T$ as in the first case 
of the conclusion of Theorem \ref{main}. 
This in particular shows that the two cases of the theorem can be 
consistent at least when we do not assume path connectedness. 

Consider the left action of the free group $F_2 =\langle a,b \rangle$ 
on its Cayley graph as an isometric action on a real tree. 
When we assume that the edges are of length 1 and we neglect the colours, 
the Cayley graph becomes a real tree, say $T$. 
Let $\theta_{ab}$ be the isometry of $T$ induced by 
the substitution $a\rightarrow b$, $b\rightarrow a$. 

Let $L\subseteq T$ be the axis of the element $ba$. 
Let $\phi$ be a loxodromic isometry of $T$ with the axis $L$ so that 
$\phi$ takes the neutral element $\Lambda\in F_2$ (the empty word) to $b$. 
In general (in particular for elements outside $L$) we define $\phi$ as follows. 
Every element $w\in F_2$ can be presented in one of the following forms: 
\begin{quote}
- $(ba)^k b^{\varepsilon} w_1$ with $k\ge 0$, $\varepsilon \in \{ 0,1\}$ 
and $(ba)^k b^{\varepsilon}$ is the longest prefix of $w$ laying on $L$; \\ 
- $(a^{-1}b^{-1})^{k} a^{\varepsilon} w_1$ with $k> 0$, $\varepsilon \in \{ 0, -1\}$ 
and $(a^{-1}b^{-1})^{k} a^{\varepsilon}$ is the longest prefix laying on $L$.  
\end{quote} 
Then in the first case $\phi$ takes $(ba)^k b^{\varepsilon} w_1$ to 
$(ba)^k b^{\varepsilon} a^{\varepsilon} b^{1-\varepsilon} w_2$, where $w_2 =\theta_{ab}(w_1 )$. 
In the second case $\phi$ takes $(a^{-1}b^{-1})^k a^{-1} w_1$ (i.e. $\varepsilon =-1$ ) 
to $(a^{-1}b^{-1})^{k} \theta_{ab} (w_1 )$ and $(a^{-1} b^{-1})^k w_1$ (i.e. $\varepsilon=0$ ) 
to $(a^{-1}b^{-1})^{k-1} a^{-1} \theta_{ab} (w_1 )$. 
It is easy to see that $\phi^2$ is the isometry of $T$ by multiplying by $ba$. 
We will apply below the easy observation that any word of even length is 
taken by $\phi$ to a word of odd length. 

Let $G$ be the subgroup of isometries of $T$ generated by $a^2$, $b^2$ and $\phi$. 
It is straightforward that $a^2 =a\phi^{-1}a^{-1}\phi$ and $b^2 = \phi a \phi^{-1}a^{-1}$. 
Thus $G$ is generated by the conjugacy class of $\phi$. \\  
{\bf Claim. } 
\begin{enumerate} 
\item {\em For every element $v\in T$ the orbit $Gv$ is contained 
in the orbit of $v$ with respect to $\langle F_2 , \theta_{ab}\rangle$. 
The median closure $c(Gv )$ contains all vertices of the Cayley graph of $F_2$. 
In particular the action is minimal. }  
\item {\em If $v\in T$ is not a vertex of the Cayley graph of $F_2$, then 
$c(Gv)$ consists of two $G$-orbits. 
For any axis $L'$ of loxodromic elements of $G$ with 
$Gv\cap L'\not=\emptyset$ the intersection $c(Gv)\cap L'$ 
does not contain $G$-orbits which are dense in $c(Gv)\cap L'$. }
\end{enumerate} 
{\em Proof of Claim.} (1) 
To see the first statement note that the $G$-orbit of the edge 
$(\Lambda ,a)$ is a subset of the orbit of $(\Lambda ,a)$ with 
respect to $\langle F_2 ,\theta_{ab} \rangle$. 
Moreover the orbit of the point $\Lambda \in T$ 
contains the points $a,b,a^{-1},b^{-1}$. 
Indeed, $a^{-1}=\phi^{-1} (\Lambda )$, $b=\phi (\Lambda )$ and 
$a = a^2 \cdot a^{-1}$, $b^{-1}=b^{-2} \cdot b$. 
This easily implies that $G\Lambda$ coincides with the Cayley 
graph of $F_2$ (which is a discrete $G$-invariant subset of $T$). 
It is easy to see that for any $v\in T\setminus L$ the median of  
$v, \phi^{-1}(v)$ and $\phi (v)$ is a vertex of the Cayley graph 
of $F_2$ (and of the axis $L$). 
If $v\in L$ then $a^{-1}=m(v, a^2 v, a^{-2}v)$ or 
$b= m(v, b^2 v, b^{-2}v)$. 
Thus all vertices of the Cayley graph of $F_2$ are included into $c(Gv)$. 

(2) By (1) the group $G$ acts without inversions 
on the set of edges of the Cayley graph of $F_2$. 
Thus when $v\in T$ is not a vertex of the Cayley graph of $F_2$ 
any edge of this graph has at most one point from $Gv$.  
On the other hand the end-points of this edge belong 
to the same orbit $G\Lambda$. 
To see that these two orbits form a median pretree 
(i.e. $c(Gv)$) it suffices to note that the median of 
a non-linear triple of the union of these orbits always 
belongs to $G\Lambda$ (the set of all points of valency $>2$). 
The rest is clear. $\Box$ 
\end{Remark}

\begin{Remark} 
It is worth noting that in the example of the previous remark 
(see Claim there) the axis $L'$ of $c(Gv)$ does not have dense 
$G_{{\bf e}}$-orbits (for appropriate ${\bf e}$) just because 
this already hods for $G$-orbits. 
On the other hand if we take $v$ to be $\Lambda$, then the axis 
of $\phi$ in  $c(Gv)$ consists of a unique orbit with respect 
to $\langle \phi \rangle$. 
It is intersting to note that for the axis $L_a$ of $a^2$ the 
situation is different. 
\begin{quote} 
{\em If ${\bf e}$ is the $+\infty$-end of $L_a$, then for any 
$v\in L_a$ the $L_a$-part of the $G_{{\bf e}}$-orbit $G_{{\bf e}}v$  
does not coincide with $c(Gv)\cap L_a$. }  
\end{quote} 
To see this statement we need few additional observations. 
Firstly note that for any element $g\in G$ there is a natural 
number $n$ and an element of the form $g'= w_1 \phi^{\varepsilon} w_2$ 
with $w_1 ,w_2 \in F_2$ of even length and $\varepsilon \in \{ -1 ,0,1 \}$ 
such that for all $a^m$ with $m>n$, $g\cdot a^m = g' \cdot a^m$.  
Indeed, $g$ can be presented in the form $w_1 \phi^{\varepsilon} g_1$ 
where $w_1\in F_2$ is of even length, $\varepsilon \in \{ -1 ,0,1 \}$ 
and $g_1\in G$ is written as a word over $\{ a,b,\phi \}$, with 
an even number of entries of $\phi$ so that each maximal $F_2$-subword 
is of even length. 
When we apply $g_1$ to a word $a^m$ with $m$ large and even, we 
obtain a word of even length having a suffix of the form $a^k$. 
Inserting a subword of the form $a^{-l} a^{l}$ if necessary we present 
the result of $g_1 \cdot a^m$ as a word $w_2 \cdot a^m$ where 
$w_2 \in F_2$ is of even length. 
This defines $g'$ as above. 

Assume that $g\in G_{{\bf e}}$. 
Then $g'$ as above belongs to $G_{{\bf e}}$ too. 
This means that when we apply $g'$ to $a^m$ with $m$ 
large enough we obtain some $a^k$. 
This only possible when $\phi^{\varepsilon} =\Lambda$ in $g'$. 
Thus $g'$ is a word from $F_2$ of even length and 
this word does not contain any entry of $b$, 
i.e. $g' = a^{2j}$ for some integer $j$. 
It is now staightforward that for any $v\in L_a$ 
the minimal distance from $v$ to another point 
of the set $G_{{\bf e}}v \cap L_a$  
is an even natural number. 
This finishes the proof. $\Box$ 
\end{Remark} 
\bigskip

\subsection{Proof of Theorem \ref{main}.}
We start with the observation that the subspace
$\bigcup\{ C_{h}:h\in X\}$ is a full subtree of $T$.
Indeed, for any $h,h'\in X$ there exist an $X$-path $h_1 ,...,h_k$ 
between $h$ and $h'$.
If $X$ consists of loxodromic elements, 
by Proposition \ref{3loxo} there is an arc in (say) 
$C_{h_{1}}\cup C_{h_{1}h_{2}}\cup ...\cup C_{h_{k}}$ joining 
$C_{h}$ and $C_{h'}$.
If $X$ consists of elliptic elements, we apply Proposition \ref{2.6}(2) 
and obtain an arc in $C_{h_1}\cup C_{h_2} \cup ... \cup C_{h_k}$ 
between $C_{h}$ and $C_{h'}$. 
Since $X$ is a conjugacy class, by minimality we see that 
$T=\bigcup\{ C_{h}:h\in X\}$.
\parskip0pt 

Denote $c(Gv )$ by $T'$. 
{\em We want to embed the $G$-pretree $T'$ into some 
special $\mathbb{R}$-tree with an isometric action of $G$.} 
If $T'$ is discrete, then it can be considered 
as a simplicial tree with an isometric action of $G$. 
Thus the pretree $T'$ can be naturally (as any simplicial tree) 
embedded into an $\mathbb{R}$-tree with an isometric action of $G$. 

It is worth noting here that when every $G_h$ is a singleton 
(in particular $h$ is elliptic), we obtain that all $C_h$ are the same 
and $G$ has a global fixed point (and in fact $T'$ and $T$ are singletons). 
From now on we assume that $T'$ is not discrete and $1< |C_h |$. 
\parskip0pt  

To prove the theorem take any $h\in X$ and an end ${\bf e}$ 
so that $C_h$ is {\bf e}-contractible. 
Assume that the first case of the conclusion of the theorem is not true.  
Choose the end ${\bf e}$ and $a_0 \in T' \cap C_h$ so that 
the $C_h$-part of the orbit $G_{{\bf e}}a_0$ is dense in $T' \cap C_h$.  
Note that now the orbit $Ga_0$ is dense in $T'$.

{\bf Claim 1.} {\em All intervals of $T'$ are dense.} 

Since $T'$ is not discrete there is an infinite sequence 
$a_1 ,...,a_k ,...\in T'$ so that each $a_{k+2}$ is between $a_1$ and $a_{k+1}$.  
Since any $C_{h'}$ with $h'\in X$, is a closed convex subset of $T$ 
we may assume that the sequence belongs to some $C_{h'}$. 
Moreover $h'$ can be chosen $h'=h$ and $a_1 ,...,a_k ,...\in G_{{\bf e}}a_0$. 

Let $g\in G_{{\bf e}}$ be loxodromic. 
We may assume that $(-\infty ,a_0 ]$ is a half-line of 
$L_g$ representing ${\bf e}$. 
If $a_1 \le_{{\bf e}} a_i$ or $a_i \le_{{\bf e}} a_1$ 
for infinitely many $i$, then applying some element of $G_{{\bf e}}$  
to $a_1$ or to appropriate $a_i$ we obtain an infinite sequence 
of vertices in some segment of $L_g$. 
If $a_1$ and $a_i$ are not ${\bf e}$-comparable for 
almost all $i$, then after replacing $a_1$ by the $\le_{{\bf e}}$-least 
element from $[a_1 ,a_2 ]\cap G_{{\bf e}}a_0$ we can apply the same 
argument to show that $L_g \cap G_{{\bf e}}a_0$ is not discrete. 

By Lemma \ref{1.12} the ordering
$L^{a_{0}}_{g} = G_{{\bf e}}a_{0}\cap L_{g}$ is dense.  
Note that in any image of $C_h$ by an element of $G_{{\bf e}}$
the corresponding part of $L^{a_0}_g$ is dense in the part 
of $T'\cap L_g$. 
Since any segment of $T' \cap C_h$ can be decomposed into two 
segments which are images of some segments from $T' \cap L_g$ 
(by ${\bf e}$-contractibility), 
we see that {\em any segment of $T'\cap C_h$ is dense}. \parskip0pt

Notice that if $h',h''\in X$ are loxodromic and $a\in L_{h'} \cap L_{h''}$  
defines an $L_{h'}$-half-line without other common elements 
with $L_{h''}$, then $a$ is the median of three non-linear 
elements from $Gv$ and thus belongs to $c(Gv)$ ($=T'$).
On the other hand if $h'$ and $h''\in X$ are elliptic 
and $a\in T^{h'} \cap T^{h''}$ defines a half-line 
without other common elements with $T^{h''}$ but with 
infinite intersection with $T^{h'}$, then $a$ is the median 
of three non-linear elements from $Gv$: two elements from $T^{h'}$, 
$T^{h''}$ and one fixed by $h'h''(h')^{-1}$. 
Thus $a\in T'$.     \parskip0pt

If now elements $a_{1}$ and $a_{2}$ belong to $T'$, say 
$a_{1}\in C_{h}$ and $a_2 \in C_{h'}$, 
then as above we find an $X$-path $h_{0},h_{1},...,h_{k}$ such that 
there is a $T$-arc in $C_{h_{0}}\cup C_{h_{1}}\cup ...\cup C_{h_{k}}$ 
joining $C_{h}$ and $C_{h'}$. 
Since $T'$ is median, the intervals of the corresponding 
$C_{h_i}$ have extremities belonging to $T'$.
This implies that the $T'$-interval $[a_1 ,a_2 ]$ is decomposed 
in $T'$ into several intervals from the corresponding $C_{h_i}$. 
In particular we now see that {\em all intervals of $T'$ are dense}. 

It is worth noting here that the argument above also shows that 
{\em the set $L^{a_0}_g$ is dense in the line $L_g \cap T'$.}  
Indeed, since $L^{a_0}_g$ is dense, any segment of $L^{a_0}_g$ 
covers $L^{a_0}_g$ by its $G_{{\bf e}}$-translations with respect 
to the action $*_g$. 
Since one of these segments is a dense subset 
of some segment of $L_g \cap T'$ we have the statement above. 

Let $T^{*}$ be the set consisting of $T'$ and all flows
of $T'$ which are induced by endless directed arcs and 
which do not have maximal elements.
We will show below that $T^{*}_0$, the set of non-terminal points of 
$T^*$ with respect to the standard betweenness relation (\cite{bow}, p.30), 
can be presented as an $\mathbb{R}$-tree where the action of $G$ is isometric.  
We start with some helpful observation. \parskip0pt 

{\bf Claim 2.} {\em Every endless directed arc $I$ from $T'$
which does not define an end of $T$, is cofinal with
an endless directed arc from some $C_{h'}\cap T'$, $h'\in X$.}
\parskip0pt 

Indeed, let $a\in I$.
Since $T$ together with the set of ends forms a complete 
$\mathbb{R}$-tree, there is $c\in T$ such 
that $I$ is cofinal with $[a,c)\cap T'$.
If $a\in C_{h'}$ and $c\in C_{h''}$, then as above we find
$h_{0},h_{1},...,h_{k}\in X$ such that there is an arc in
$C_{h_{0}}\cup C_{h_{1}}\cup ...\cup C_{h_{k}}$ joining $C_{h'}$
and $C_{h''}$.
We see that $[a,c)$ is decomposed into finitely many 
intervals from the corresponding $C_{h_i}$.
The last interval (which is of the form $[a',c)$ with 
$a'\in T'$) is cofinal with $I$.

We now extend the betweenness relation to $T^{*}$ as in \cite{bow}
(p.30): for points $x,y$ and a flow $r$ we say $B(y;x,r)$ if
$(x,y) \in r$.
Then we can define $B(r;x,y)$ as the case when
there is no point $z$ with $B(z;x,r)\wedge B(z;y,r)$.
For flows $r,r'$ and a point $x$ we say $B(x;r,r')$ if for
any point $y\not= x$, $(y,x)\in r$ or $(y,x)\in r'$.
We also put $B(r;x,r')$ if $B(r;x,y)$ for some $y\in T'$ with
$(x,y)\in r'$.
If $p,q,r$ are flows then $B(p;q,r)$ means that there is no
point $z$ with $B(z;p,q)\wedge B(z;p,r)$. \parskip0pt 

Since we consider flows without maximal elements 
this definition implies that any two distinct flows 
from $T^*\setminus T'$ are separated by an element from $T'$. 
Thus for {\em a non-linear triple
$u,v,w\in T^{*}$ the median $m(u,v,w)$ can be presented
as $m(x,y,z)$ for some non-linear $x,y,z\in T'$;
thus $m(u,v,w)\in T'$}.\parskip0pt 

By Claim 2 we know that {\em all flows of $T^* \setminus T'$ 
are induced by endless directed arcs of} $C_{h'} \cap T'$, 
$h'\in X$. \parskip0pt 

It is also worth noting that the definition of the betweenness 
relation on $T^*$ implies that a flow $r\in T^* \setminus T'$ 
{\em lies on a line $L'\subset T'$ if and only if $r$ forms 
a linear triple with any two points of} $L'$.\parskip0pt 
 
As any flow $r\in T^{*}\setminus T'$ is induced by any of its
linear orderings without upper bounds (Lemma \ref{flowcut}), 
if $r$ belongs to a $T^{*}_0$-line $L^{*}$, then $r$ is defined 
by a class of the corresponding Dedekind cut of $L'=T'\cap L^{*}$.
Moreover any cut $L'= C^{-}\cup C^{+}$, which does not
define an element of $T'$, defines a flow from $T^{*}_0$.
Indeed suppose for a contradiction that $r$ is the flow 
defined by $C^{-}$ and $c\in T'$ is a maximal element of $r$.
Then for any $c'\in C^{-}$ and $c''\in C^{+}$,
the interval $(c',c)$ in $T'$ consists of elements 
of $L'$ $r$-greater than $c'$ (by axioms of flows and 
maximality of $c$) and is contained in $[c',c'']$ 
(see Lemma \ref{flowcut}). 
Find the median $m(c,c',c'')\in T'$. 
If $c\not= m(c,c',c'')$, then by density the interval $(m(c,c',c''),c)$ 
is infinite, a contradiction with $(c',c) \subseteq [c',c'']$. 
Thus, $c\in [c',c'']$,
contradicting the assumption that the cut does not 
define an element of $T'$.
As a result we obtain that {\em the $T^{*}_0$-line
$L^{*}$ is the Dedekind completion of $L'$ and can be 
identified with} $\mathbb{R}$. \parskip0pt

Since $T'$ is dense, $T'$ {\em is dense in} $T^*_0$.\parskip0pt 
 
The action of $G$ on $T'$ uniquely defines an action on $T^{*}_0$.
We want to show that this action is an isometric action on
an $\mathbb{R}$-tree.
Let us start with the following claim. \parskip0pt 

{\bf Claim 3.} {\em If $g'\in G$ is not loxodromic in $T$, 
then $T'$ contains a point fixed by $g'$ or a segment 
which is inversed by $g'$.} \parskip0pt 

Indeed, assume that $g'$ does not fix any point of the orbit $Gv$.  
Find $b\in Gv$ and $a\in T$ so that $[b,g'(b)]\cap T^{g'} =a$ 
(see Lemma \ref{lem_nnest} and Proposition \ref{2.6}). 
If the interval $(a,g'g'(b))$ does not meet $(a,b)\cup (a,g'(b))$ 
then $a$ is the median of $b,g'(b)$ and $g'g'(b)\in Gv$, i.e. 
$a\in T'=c(Gv)$. 
If $(a,g'(b))\cap (a,g'g'(b))\not=\emptyset$, 
then $(a,b)\cap (a,g'(b))\not=\emptyset$, a contradiction. 
In the case $(a,b)\cap (a,g' g' (b))\not=\emptyset$  
we see that when $g'g'(b)\not= b$, the median $c$ of $b,g'(b)$ and $g'g'(b)$ 
belongs to $(a,b)\cap (a,g' g' (b))\cap c(Gv )$. 
By non-nesting, $g' g' (c)=c$. 
Thus $[c,g' (c)]$ is inversed by $g'$. \parskip0pt

{\bf Claim 4.} {\em An element $g'\in G$ fixes a point in $T$ 
if and only if it fixes a point in $T^{*}_0$.} \parskip0pt

If $g'$ does not fix a point in $T$, then there is a line 
$L\subseteq T$ which is the axis of $g'$.
Since $C_h$ is ${\bf e}$-contractible and every segment of $T$ can 
be presented as the union of finitely many segments from $C_{h'}$, 
$h'\in X$, with common extremities belonging to $T'$, 
the set $L\cap T'$ is not empty and thus is cofinal with $L$. 
By the definition of $T^{*}_0$ and the corresponding betweenness 
relation, $L\cap T'$ is cofinal with some $T^{*}_0$-line 
$L^{*}$.
It is clear, that $L^{*}$ is the axis of $g'$ in $T^{*}_0$.
Thus $g'$ is loxodromic in $T^{*}_0$. 
\parskip0pt 

To see the converse we apply Claim 3.
By this claim we only have to consider the case when $g'$ 
fixes a point in $T$, does not fix any point of $T'$ and 
inverses a segment, say $[c,g'(c)]$, in $T'$. 
Let $I=\{ x\in [c,g'(c)]\cap T' :B(x;c,g'(x))\}$.
Since $T'$ is dense, $I$ is an endless directed arc. 
It is straightforward, that $I$ and $g'(I)$ define the same 
flow which is fixed by $g'$.     
\parskip0pt 

{\bf Claim 5.} {\em The action of the group $G$ on 
$T^{*}_0$ is non-nesting.}\parskip0pt 

Indeed, if $g'\in G$ is loxodromic in $T$, then it is loxodromic
in $T^{*}_0$. 
Lemma \ref{loxo} implies that if $g'$ maps a segment of $T^{*}_0$ 
properly into itself, then such a segment can be chosen in the 
axis $L^{*}_{g'}$ of $g'$. 
The latter is impossible, because $g'$ has a non-nesting action 
on $L^{*}_{g'} \cap T'$ with a cofinal orbit and $L^{*}_{g'}\cap T'$ 
is dense in $L^{*}$.   

If $g'$ fixes a point $a\in T^{*}_0$ and maps a segment $[b,b']$ 
properly into itself then $a$ and $[b,b']$ can be chosen so 
that $a\in [b,b']$. 
This can be shown by straightforward arguments depending 
on the place where the median of $a$ and the extremities of 
the segment lie. 
If it happens that $a=b'$, then as $T'$ is dense in $[a,b]$, 
we can arrange that $b\in T'$.
Since the action of $G$ on $T'$ is non-nesting we see that 
$a\not\in T'$.
Using Claims 3 and 4 find a segment $[c,g'(c)]$ with 
$(g')^{2}(c)=c\in c(Gv )$ and $a\in [c,g'(c)]$.
Since $a\not\in T'$, $a$ is not the median of $b,c,g'(c)$. 
Thus one of the intervals $(a,m(a,b,c))$ or $(a,m(a,b,g'(c))$ 
is non-empty and contains a point from $T'$. 
Replacing $c$ by this point if necessary,
we may assume that $c\in [a,b]$. 
Then $g'$ maps $[g'(c),b]$ properly into itself, contradicting 
non-nesting on $T'$. \parskip0pt 

Consider the case when neither $b$ nor $b'$ are fixed by $g'$.
If $g'(b)\in [a,b]$ or $g'(b')\in [a,b']$ then apply 
the argument of the previous paragraph.      
Assume $g'(b)\in [a,b']$. 
Then $g'(b')\in [a,b]$ and $(g')^2$ maps the segment $[a,b']$ 
properly into itself.
As we already know this contradicts the assumption that the 
action of $G$ on $T'$ is non-nesting. 
This finishes the proof of the claim.

Let $g$, ${{\bf e}}$ and $h$ be as in the begining of the proof. 
Let $s\in G$. 
As we already know, the line of $T^{*}_0$ containing a copy
$s^{-1}(L^{a_{0}}_{g})$ can be identified with
the Dedekind completion of the $T'$-line  
$s^{-1}(T' \cap L_{g})$.
Since $L^{a_0}_{g}$ is dense in $T' \cap L_g$, this 
completion coincides with the Dedekind completion  
$s^{-1}(L^{a_0}_{g})^{*}$ of $s^{-1}(L^{a_0}_g )$. 
It is clear that $s^{-1}(L^{a_0}_g )^{*}$ is the axis 
$L^*_{s^{-1}gs}$ of $s^{-1}gs$ in $T^*_0$. 
Note that $L^*_{s^{-1}gs}= s^{-1}(L^*_g )$. \parskip0pt

It is easy to see that $G_{{\bf e}}$ coincides with the stabiliser 
in the $G$-space $T^*_0$ of the end represented by a half-line of 
the $g$-axis $L^*_g$ containing $(-\infty, a_0 ]\cap T'$. 
Since every $g'\in G_{{\bf e}}$ acts on $L_g \subset T$ as 
a translation (up to topological equivalence with $\mathbb{R}$) 
under the action $*_{g}$ of $G_{{\bf e}}$ (as in Section 2.2),  
so it does on the axis $L^{*}_g$ (by non-nesting). 
It is clear that all elements of $G_{({\bf e})}$ 
fix $L^{*}_g$ pointwise.     \parskip0pt 
          
The line $L^{*}_g$ can be identified 
with $\mathbb{R}$ so that the group $L^{a_{0}}_{g}$ 
(under the structure of the Archimedean group 
$G_{{\bf e}}/G_{({\bf e})}$ defined in Section 2.2) acts on
$L^{*}_{g}$ by translations defined by real numbers.
We assume that $a_{0}$ corresponds to $0$.
Let $d$ denote the metric obtained on $L^*_{g}$ in this way.
\parskip0pt

Now for any $s\in G_{{\bf e}}$ consider $s^{-1}(L^*_{g})$ under
the metric $d^{s}$ induced by $d$ and the map $s^{-1}$ so that 
$s^{-1}$ is an isometry. \parskip0pt 

{\bf Claim 6.} {\em For any $s,s'\in G_{{\bf e}}$ the metrics $d^{s}$ 
and $d^{s'}$ agree on the common half-line of the corresponding lines.
In particular, if $s^{-1}(L^*_g ) =(s')^{-1}(L^*_g )$, 
then $d^s =d^{s'}$. } \parskip0pt

Suppose $x_{1},x_{2}\in s^{-1}(L^*_{g})\cap (s')^{-1}(L^*_{g})$,
$x_{1} < x_{2}$ and $d^{s'}(x_{1},x_{2}) < d^{s}(x_{1},x_{2})$
(the case $d^{s}(x_{1},x_{2}) < d^{s'}(x_{1},x_{2})$ is similar).
Since $(s')^{-1}(L^{a_{0}}_{g})$ is dense in $(s')^{-1}(L^*_{g})$, 
we can find $x'_{1}\le x''_{1} < x''_{2}\le x'_{2}\in s^{-1}(L^*_{g})\cap (s')^{-1}(L^*_{g})$
with  $d^{s'}(x'_{1},x'_{2}) < d^{s}(x''_{1},x''_{2})$ 
where $x'_{1},x'_{2} \in (s')^{-1}(L^{a_{0}}_{g})$ and 
$x''_{1},x''_{2} \in s^{-1}(L^{a_{0}}_{g})$. 
Then the inequality $s'(x'_{2})-s'(x'_{1})<s(x''_{2}) - s(x''_{1})$
holds in the group $L^{a_{0}}_{g}$. 
Applying a translation from $L^{a_{0}}_{g}$ we can take
$s(x''_{1})$ to $s'(x'_{1})$; then $s(x''_{2})$ must go to an
element of $L^{a_{0}}_{g}$ greater that
$s'(x'_{2})$, a contradiction with non-nesting.

Extend the metric $d$ to the space 
$T^*_{{\bf e}} = \bigcup \{ s^{-1}(L^*_{g}): s\in G_{{\bf e}}\}$ 
as follows. 
If $a\in s^{-1}(L^*_{g})$ and $b\in (s')^{-1}(L^*_{g})$ 
both belong to one of the lines $s^{-1}(L^*_{g})$ or 
$(s')^{-1}(L^*_{g})$ then the distance between them is 
defined by the distance in the corresponding axis. 
Otherwise find the arc $[a,c]\cup [c,b]$ where 
$(-\infty, c]=s^{-1}(L^*_{g})\cap (s')^{-1}(L^*_{g})$.
Now define the distance between $a$ and $b$ as the sum 
of distances $d^{s}(a,c)$ and $d^{s'}(c,b)$. 

To prove that $d$ on $T^*_{{\bf e}}$ is invariant 
under the $G_{{\bf e}}$-action it suffices to show that 
if $s\in G_{{\bf e}}$ maps $s^{-1}_{1}(L^*_{g})$
onto $s^{-1}_{2}(L^*_{g})$, then $s$ maps the metric
$d^{s_{1}}$ onto $d^{s_{2}}$.
The latter condition can be verified as follows. 
By Claim 6 the metrics $d^{s_1}$ and $d^{s_2 s}$ 
coincide on $s^{-1}_1 (L^*_g )$.
The latter one is the image of $d^{s_2}$ under $s^{-1}$
by the definition. 

Let $C^*_h$ be the the characteristic set of $h$ 
in $T^*_0$ ($h$ is as above). 

{\bf Claim 7.} $C^*_h \subseteq T^*_{{\bf e}}$. 

If $h$ is loxodromic in $T$, then by Claim 4, $h$ 
is loxodromic in $T^*_0$ and thus the axis $C^*_h$ 
is cofinal with $C_h$. 
Since $C_h$ is ${\bf e}$-contractible, the set 
$C^*_h$ is ${\bf e}$-contractible too. 
As a result $C^*_h \subset T^*_{{\bf e}}$. 

When $h$ is elliptic in $T$, by Claim 4 it is elliptic in $T^*_0$. 
If $w\in (T')^h$, then 
$w\in T^*_{{\bf e}}$ by the assumptions of the theorem.  
If $h$ fixes $r \in T^*_0 \setminus T'$, then 
as we already know $r$ {\em cannot be the median of three 
non-linear points in $T^*_{0}$.}  
Let $I$ be an endless directed arc in $T'$ inducing the flow $r$. 
By Lemma \ref{flowcut} there are only two possibilities: there is 
a cofinal arc in $I$ which is fixed pointwise by $h$ or 
$r$ is between $u$ and $h(u)$ for any $u\in I$. 
In the first case find $w\in T^h$ determined by $I$ as the end-point of 
the corresponding interval. 
As $T^h \subseteq T_{{\bf e}}$, there is $s\in G_{{\bf e}}$ with $w\in s(L_g )$. 
Since $s(L_g )$ has a non-empty intersection with $T'$, $w\not\in T'$ and 
$r$ is not the median of a non-linear triple, 
a cofinal arc of $I$ is contained in $s(L_{g})$. 
Thus the point $r$ belongs to $s(L^*_g )\subseteq T^*_{{\bf e}}$ too.  
In the second case there is a point 
$w\in T^h$ which belongs to some $s(L_g )$ with $s\in G_{{\bf e}}$ 
and which is between $u$ and $h(u)$ for any $u\in I$. 
If a cofinal part of $I$ belongs to $s(L_g )$, then $r\in s(L^*_g )$ too. 
Otherwise since $s(L_g )$ has a non-empty intersection with $T'$ 
there is $u'\in T'$ which is the projection of any $u\in I$ onto $s(L_g )$ 
(it is the median of a non-linear triple from $T'$ and thus belong to $T'$). 
Thus $u'$ is between any $u\in I$ and $h(u)$. 
This contradicts to the assumption that $r$ is a flow defined by the endless 
directed arc $I$ and also is defined by the endless directed arc $h(I)$. 
As a result we conclude that {\em if $h$ is elliptic the characteristic 
set $C^*_h$ is contained in $T^*_{{\bf e}}$.} 

We now see that there is a metric $d$ on $C^*_h$ which is 
$G_{{\bf e}}$-invariant, i.e. satisfies the property that 
if $c_1 ,c_2 \in C^*_h$, $s\in G_{{\bf e}}$ and $s(c_1 ), s(c_2 )\in C^*_h$, 
then $d(c_1 ,c_2 )=d(s(c_1 ),s(c_2 ))$.  
Moreover by inspection we also see that such a metric can 
be obtained from any $G_{{\bf e}}$-invariant metric on $L^*_g$ 
with respect to the action $*_g$ (see Remark \ref{more}). 

Now for any $s\in G$ consider $s^{-1}(C^*_{h})$ under
the metric $d^{s}$ induced by $d$ and the map $s^{-1}$ so that 
$s^{-1}$ is an isometry. \parskip0pt 

{\bf Claim 8.} {\em For any $s,s'\in G$ the metrics 
$d^{s}$ and $d^{s'}$ agree on $s^{-1}(C^*_{h}) \cap (s')^{-1}(C^*_h )$.  
In particular, if $s^{-1}(C^*_h ) =(s')^{-1}(C^*_h )$, 
then $d^s =d^{s'}$. } \parskip0pt

Suppose
$x_{1},x_{2}\in s^{-1}(C^*_h )\cap (s')^{-1}(C^{*}_{h})$,
and $d^{s}(x_{1},x_{2}) < d^{s'}(x_{1},x_{2})$
(the case $d^{s'}(x_{1},x_{2}) < d^{s}(x_{1},x_{2})$ is similar).
As we already know any segment of $C^*_h$ can be decomposed 
into a sum of two subsegments so that each of them can be 
taken to $L^*_g$ by an element of $G_{{\bf e}}$. 
We may suppose that $[x_1 ,x_2 ]$ can be taken to $s^{-1}(L^*_g )$ 
by an element of $G_{s^{-1}({\bf e})}$ and can be taken to 
$(s')^{-1}(L^*_g )$ by an element of $G_{(s')^{-1}({\bf e})}$. 

Find
$x'_{1}, x'_{2}, x''_{1}, x''_{2} \in s^{-1}(C^*_h )\cap (s')^{-1}(C^{*}_{h}) \cap T'$
with  $[x''_1 ,x''_2 ]\subset [x'_1 ,x'_2 ]$ (where $x''_2$ is between $x''_1$ and $x'_2$) 
and $d^{s}(x'_{1},x'_{2}) < d^{s'}(x''_{1},x''_{2})$. 
Since $L^{a_{0}}_{g}$ is dense in $L^{*}_{g}$, the points 
$x'_{1}, x'_{2}, x''_{1}, x''_{2}$ can be chosen so that for appropriate 
$s_1 \in G_{s^{-1}({\bf e})}$ and $s'_1 \in G_{(s')^{-1}({\bf e})}$ 
the inequality $ss_1 (x'_{2})-ss_1 (x'_{1})<s's'_1 (x''_{2}) - s's'_1 (x''_{1})$
holds in the group $L^{a_{0}}_{g}$.
Applying a translation from $L^{a_{0}}_{g}$ we can take
$s's'_1 (x''_{1})$ to $ss_1 (x'_{1})$; then $s's'_1 (x''_{2})$ must 
go to an element of $L^{a_{0}}_{g}$ greater that
$ss_1 (x'_{2})$, a contradiction with non-nesting.

Since any segment of $T^*_0$ is decomposed into finitely many 
segments from appropriate $C^*_h$, $h\in X$, we can extend all metrics 
$d^s$, $s\in G$, to a metric $d$ on the space $T^*$.
It is straightforward that $d$ defines an $\mathbb{R}$-tree. 
By Claim 8 to prove that $d$ is invariant under the $G$-action 
it suffices to show that if $s$ maps $s^{-1}_{1}(C^{*}_{h})$
onto $s^{-1}_{2}(C^{*}_{h})$, then $s$ maps the metric
$d^{s_{1}}$ onto $d^{s_{2}}$.
The latter condition can be verified as follows. 
By Claim 8 the metrics $d^{s_1}$ and $d^{s_2 s}$ coincide on 
$s^{-1}_1 (C^{*}_h )$.
The latter one is the image of $d^{s_2}$ under $s^{-1}$
by the definition.  
 $\Box$

\end{document}